\documentclass[a4paper,12pt]{article}
\usepackage[english]{babel}
\usepackage[latin1]{inputenc}
\usepackage[T1]{fontenc}
\usepackage{amsmath}
\usepackage{amsfonts}
\usepackage{amsthm}
\usepackage{amssymb}
\usepackage{dsfont}
\usepackage{epsfig}
\psfigdriver{dvips}
\newtheorem{theorem}{Theorem}[section]
\newtheorem{proposition}[theorem]{Proposition}
\newtheorem{remark}[theorem]{Remark}

\newtheorem{corollary}[theorem]{Corollary}

\newcommand{\R}{\mbr} 
\newcommand{\absj}[1]{\left\lvert #1 \right\rvert} 
\newcommand{\1}{\mathds{1} } 
\newcommand{\Proba}{\mbp}
\newcommand{\E}{\mathbb{E}} 
\DeclareMathOperator{\card}{Card} 
\DeclareMathOperator{\var}{Var} 

\newcommand{\telque}{\mbox{ s.t. }}

\newcommand{\mbe}{\mathbb{E}}
\newcommand{\mbr}{\mathbb{R}}
\newcommand{\mbp}{\mathbb{P}}

\newcommand{\mbt}{\mathbb{T}}
\newcommand{\mbz}{\mathbb{Z}}

\newcommand{\mtc}{\mathcal}
\newcommand{\mbf}{\mathbf}

\newcommand{\wt}[1]{{\widetilde{#1}}}

\newcommand{\ol}[1]{\overline{#1}}

\newcommand{\e}[1]{\mbe\brac{#1}}
\newcommand{\ee}[2]{\mbe_{#1}\brac{#2}}
\newcommand{\prob}[1]{\mbp\brac{#1}}
\newcommand{\probb}[2]{\mbp_{#1}\brac{#2}}
\newcommand{\paren}[1]{\left(#1\right)}
\newcommand{\brac}[1]{\left[#1\right]}

\newcommand{\norm}[1]{\left\|#1\right\|}
\newcommand{\set}[1]{\left\{#1\right\}}
\newcommand{\abs}[1]{\left| #1 \right|}

\newcommand{\defeq}{:=}

\newcommand{\al}{\alpha}

\newcommand{\Om}{\Omega}

\newcommand{\gru}{W} 

\newcommand{\bY}{{\mbf{Y}}}
\newcommand{\bA}{{\mbf{A}}}
\newcommand{\bG}{{\mbf{G}}}

\newcommand{\cB}{{\mtc{B}}}

\newcommand{\moy}[1]{\overline{#1}} 
\newcommand{\moyW}[2]{\moy{#1}_{\left[ #2 \right]}} 

\newcommand{\latin}[1]{\textit{#1}}

\DeclareMathOperator{\supp}{supp} 

\begin{document}
\title{Resampling-based confidence regions and multiple tests for a correlated random vector}
\author{
Sylvain Arlot\\
Univ Paris-Sud, Laboratoire de Math\'ematiques d'Orsay, \\
Orsay Cedex, F-91405; CNRS, Orsay cedex, F-91405 \\
\texttt{sylvain.arlot@math.u-psud.fr}
\\
INRIA Futurs, Projet Select
\and
Gilles Blanchard\\
Fraunhofer FIRST.IDA, Berlin, Germany,\\
\texttt{blanchar@first.fraunhofer.de}
\and
\'Etienne Roquain\\
INRA Jouy-en-Josas, unit\'e MIG, \\
78 352 Jouy-en-Josas Cedex, France,
\\
\texttt{etienne.roquain@jouy.inra.fr}
}
\maketitle

\begin{abstract}
We derive non-asymptotic confidence regions for the mean of a random vector
whose coordinates have an unknown dependence structure. The random
vector is supposed to be either Gaussian or to have a symmetric bounded distribution, and we
observe $n$ i.i.d copies of it. The confidence regions are built
using a data-dependent threshold based on a weighted bootstrap
procedure. We consider two approaches, the first based on a
concentration approach and the second on a direct boostrapped 
quantile approach. The first one allows to deal with 
a very large class of resampling weights while our results for
the second are restricted to Rademacher weights. However, 
the second method seems more accurate in practice. Our results
are motivated by multiple testing problems, and we show on
simulations that our procedures are better than the
Bonferroni procedure (union bound) as soon as the observed vector has
sufficiently correlated coordinates.
\end{abstract}



\section{Introduction}

In this work, we assume that we observe a sample $\bY \defeq (\bY^{1},\dots,\bY^{n})$
of $n\geq 2$  
i.i.d. observations of an integrable 
random vector $\bY^i \in \R^K$
with a dimension $K$ possibly much greater than $n$. Let $\mu\in\mbr^K$ denote the common
mean of the $\bY^{i}$\,; 
our main goal is to find a non-asymptotic 
$(1-\alpha)$-confidence region for $\mu$\,, 
of the form:
\begin{equation}
\left\{ x\in  {\R}^{K} \telque \phi \paren{\moy{\bY}-x }\leq
  t_{\alpha}(\bY) \right\} \, ,\label{Equ_CR1} 
\end{equation}
where $\phi:{\R}^{K}\rightarrow {\R}$ is a measurable 
function fixed in advance by
the user (measuring a kind of distance), $\alpha\in (0,1)$, $t_{\alpha} :
\paren{{\R}^{K}}^{n} \rightarrow {\R}$ is a measurable 
data-dependent threshold,
and $\moy{\bY}=\frac{1}{n}\sum_{i=1}^{n} \bY^{i}$
is the empirical mean of the sample $\bY$.  

The form  
of the confidence region \eqref{Equ_CR1} is
motivated by the following multiple testing problem: if we want to
test simultaneously for all $1 \leq k \leq K$ the hypotheses $H_{0,k}=\{\mu_{k}\leq 0\}$
against $H_{1,k}=\{\mu_{k}> 0\}$, we propose to reject the $H_{0,k}$
corresponding to 
\[ \{1 \leq k \leq K \telque \overline{\bY}_{k} > t_{\alpha}(\bY)\} \, . \]

The error of this multiple testing procedure can be measured by the family-wise 
error rate defined by the probability that at least one hypothesis is wrongly
rejected. Here, this error will be strongly (i.e. for any value of $\mu$) 
controlled by $\alpha$ as soon as the confidence region
(\ref{Equ_CR1}) for $\mu$ 
with $\phi=\sup(\cdot)$  is of level at least
$1-\alpha$. Indeed, for all $\mu$, 
\begin{align*}
\Proba \paren{ \exists k \telque \moy{\bY}_{k}>t_{\alpha}(\bY)\mbox{ and
} \mu_{k}\leq 0 } &\leq  \Proba \paren{ \exists k \telque
\overline{\bY}_{k}-\mu_{k}>t_{\alpha}(\bY) } \\ 
&= \Proba \paren{ \sup_k \left\{ \moy{\bY}_k-\mu_k \right\} >t_{\alpha}(\bY) } \enspace .
\end{align*}
The same reasoning with $\phi=\sup\absj{\cdot}$
allows us to test $H_{0,k}=\{\mu_{k}=0\}$  
against $H_{1,k}=\{\mu_{k}\neq 0\}$, by choosing the rejection set 
$\{1\leq k\leq K \telque \: \absj{ \moy{\bY}_k } > t_{\alpha}(\bY)\}$.   

While this goal is statistical in motivation, to tackle it we want to follow
a point of view inspired from learning theory, in the following
sense: first, we want a non-asymptotical result valid for any fixed $K$ and $n$,
and secondly, we want to make no assumptions on the dependency structure of the
coordinates of $\bY^i$ (although we will consider some general assumptions
over the distribution of $\bY$, for example that it is Gaussian).

The ideal threshold $t_{\alpha}$ in \eqref{Equ_CR1} is obviously
the $1-\alpha$ quantile of the distribution of
$\phi\paren{\moy{\bY}-\mu}$. However, this quantity depends on the
unknown dependency structure of the coordinates of $\bY^i$ and is
therefore itself unknown.

We propose here to approach $t_{\alpha}$ by some resampling scheme: the
heuristics of the resampling method (introduced by Efron \cite{Efr:1979}) is that the distribution of
$\moy{\bY}-\mu$ is ``close'' to the one of  
\[ \moyW{\bY}{W - \moy{W}} :=\frac{1}{n}\sum_{i=1}^{n} (W_{i}-\moy{W})
\bY^{i} = \frac{1}{n}\sum_{i=1}^{n} W_{i}(\bY^{i}-\moy{\bY})=\moyW{\paren{\bY - \moy{\bY}}}{W} 
\, ,\]
conditionally to $\bY$, where $(W_{i})_{1 \leq i \leq n}$ are real
                                random variables independent 
of $\bY$ called the \emph{resampling weights}, and $\moy{W} = n^{-1} \sum_{i=1}^n W_i$\,.
We emphasize that the family $(W_i)_{1\leq i \leq n}$ itself {\em need not be independent.}

Following this idea, we propose two different approaches to obtain non-asymptotic confidence regions in this paper:
\begin{enumerate}  
\item The expectations of $\phi\paren{\moy{\bY}-\mu}$ and
  $\phi\paren{\moyW{\bY}{W - \moy{W}}}$ can be precisely compared, and the
  processes $\phi\paren{\moy{\bY} - \mu}$ and $\E \left[
    \phi\paren{\moyW{\bY}{W - \moy{W}}} \big| \bY \right]$ concentrate
  well around their expectations. 
\item The $1-\alpha$ quantile of the distribution of
  $\phi\paren{\moyW{\bY}{W - \moy{W}}}$ conditionally to $\bY$ is
  close to the one of $\phi\paren{\moy{\bY} - \mu}$\,. 
\end{enumerate}
Method 1 above is closely related to the Rademacher complexity
approach in learning theory, and our results in this direction are
heavily inspired by the work of Fromont \cite{Fro:2004}, who studies
general resampling schemes in a learning theoretical setting. It may also 
be seen as a generalization of cross-validation methods. 
For method 2, we will restrict ourselves specifically to Rademacher weights
in our analysis, because we use a symmetrization trick.
Although this kind of method is not new in the
resampling literature, to our knowledge our result is the
first to provide a non-asymptotic analysis based on 
empirical resampled quantiles.

Let us now define a few notations that will be useful throughout this paper.
\begin{itemize}
\item Vectors, such as data vectors $\bY^i=(\bY^i_k)_{1 \leq k \leq
  K}$, will always be column vectors. Thus, $\bY$ is a $K \times n$
  data matrix.  
\item If $\mu \in \R^K$, $\bY - \mu$ is the matrix obtained by
  subtracting $\mu$ to each (column) vector of $\bY$.  
If $c \in \R$ and $W \in \R^n$, $W-c = (W_i - c)_{1 \leq i \leq n} \in \R^n$.
\item $\overline{\Phi}$ is the standard Gaussian upper tail function.
\end{itemize}
Several properties may be assumed for the function $\phi : \R^K \rightarrow \R$:
\begin{itemize}
\item Subadditivity: $\forall x,x^{\prime} \in \R^K, \quad \phi
  \paren{ x + x^{\prime} } \leq \phi(x) + \phi \paren{ x^{\prime}
  }$\,. 
\item Positive-homogeneity: $\forall x \in \R^K, \, \forall \lambda
  \in \R_+, \quad \phi \paren{ \lambda x  }  = \lambda \phi(x)$\,. 
\item Bounded by the $p$-norm, $p\in [1,\infty]$: $\forall x \in
  \R^K$, $\absj{\phi \paren{ x  }} \leq \norm{x}_p$, where $\norm{x}_p$ is equal to $(\sum_{k=1}^{K}\absj{x_{k}}^{p})^{1/p}$ if $p <
  \infty$ and $\max_{k}\{\absj{x_k}\}$ otherwise.  
\end{itemize}

\noindent Finally, different assumptions on the
generating distribution of $\bY$ can be made:
\begin{itemize} 
\item[(GA)] The Gaussian assumption: the $\bY^{i}$ are Gaussian vectors
\item[(SA)] The symmetric assumption: the $\bY^{i}$ are symmetric
  with respect to $\mu$ i.e. $\bY^{i}-\mu \sim \mu-\bY^{i}$\,.
\item[(BA)]\hspace{-2mm}($p,M$) The bounded assumption: $\norm{\bY^{i} - \mu}_{p}\leq M$ a.s.
\end{itemize}
In this paper, our primary focus is on the Gaussian framework (GA), because the
corresponding results will be more accurate.

The paper is organized as follows: Section~2 deals with the
concentration method with general weights. In Section~3, we propose an approach based on
resampling quantiles, with Rademacher weights. 
We illustrate our methods in Section~4 with a simulation study. The proofs of our
results are given in Section~5. 

\section{Confidence region using concentration} 
 
In this section, we consider a general $\R^n$-valued \emph{resampling weight vector}
$W$\,, satisfying the following properties: $W$
is independent of $\bY$, for all $i \in \{1, \ldots, n\}$ $\e{W_i^2}<\infty$\,, the $(W_i)_{1\leq i\leq n}$ have an exchangeable distribution (\latin{i.e.} invariant under every permutation of the indices) and
the coordinates of $W$ are not a.s. equal, \latin{i.e.} $\E \absj{W_1-\moy{W}}>0$.
Several examples of resampling weight vectors are given in Section
\ref{sec:conc_poids}, where we also tackle the question of choosing a
resampling.
 
Four constants that depend only on the distribution of $W$ appear in the
results below (the fourth one is defined only for a particular class
of weights). They are defined as follows and computed
for classical resamplings in Tab.~\ref{tab:poids_ctes}:
\begin{align} 
A_W &
\defeq \E{\absj{W_1 - \moy{W}}} \label{eq:AW} \\ 
B_W &\defeq \E \left[ \paren{\frac{1}{n}\sum_{i=1}^{n} 
    \paren{W_{i}-\moy{W}} ^{2} }^{\frac{1}{2}} \right]\, \label{eq:BW} \\  
C_W & \defeq   \paren{\frac{n}{n-1} \E \left[ \paren{ W_{1}-\moy{W} }^2 \right] }^{\frac{1}{2}}\\ 
D_W &\defeq a + \E \absj{\moy{W} - x_0} 
\quad \mbox{if } \forall i, \, \absj{W_i-x_0} = a \mbox{ a.s. (with $a>0,x_0\in\R$)}\,. 
 \label{eq:DW}  
\end{align}

Note that under our assumptions, these quantities are
positive. Moerover, if the weights are i.i.d., $C_W =
\var(W_1)^{\frac{1}{2}}$. We can now state the main result of this
section:  

\begin{theorem}\label{th_mainconc} Fix $\alpha\in (0,1)$ and
  $p\in[1,\infty]$. Let $\phi:\R^{K}\rightarrow \R$ be any function
  subadditive, positive-homogeneous and bounded by the $p$-norm, and
  let $W$ be a resampling weight vector.
\begin{enumerate} 
\item If $\bY$ satisfies (GA), then
\begin{equation} 
  \phi \paren{\moy{\bY}-\mu} < \frac{ \E \left[ \phi
        \paren{\moyW{\bY}{W-\moy{W}}} \big| \bY \right] }{B_W} +
    \norm{\sigma}_{p} { \overline{\Phi}}^{-1}(\alpha/2) \left[
      \frac{C_W}{nB_W} + \frac{1}{\sqrt{n}}\right] 
  \label{Equ_concfinalgauss1}  
\end{equation} 
holds with probability at least $1-\alpha$, 
where $\sigma$ is the vector $[\var^{1/2}(\bY^{1}_{k})]_k$. The same
bound holds for the lower deviations, i.e. with inequality \eqref{Equ_concfinalgauss1} reversed and the additive term replaced by its opposite.
\item If $\bY$ satisfies (BA)($p,M$) and (SA), then 
\begin{equation}  
\phi \paren{\moy{\bY}-\mu} < \frac{ \E \left[ \phi
      \paren{\moyW{\bY}{W-\moy{W}}} \big| \bY \right] }{A_W} +
  \frac{2M}{\sqrt{n}}\sqrt{\log(1/\alpha)} 
\nonumber
\end{equation} 
holds with probability at least $1-\alpha$\,.
If moreover the weights satisfy the assumption of \eqref{eq:DW}, then 
\begin{equation} 
\phi \paren{\moy{\bY}-\mu} > \frac{ \E \left[ \phi
      \paren{\moyW{\bY}{W-\moy{W}}} \big| \bY \right] }{D_W} -
  \frac{M}{\sqrt{n}}\sqrt{1+\frac{A_{W}^{2}}{D_{W}^{2}}}\sqrt{2\log(1/\alpha)}
\nonumber 
\end{equation} 
holds with probability at least $1-\alpha$\,.
\end{enumerate} 
\end{theorem} 
 
If there exists a deterministic threshold $t_{\alpha}$ such that
$\Proba(\phi \paren{\moy{\bY}-\mu} > t_{\alpha})\leq \alpha$,
the following corollary establishes that we can combine the above concentration
threshold with $t_{\alpha}$ to get a new threshold almost better than
both.  

\begin{corollary}
\label{cor_mainconc} Fix $\alpha,\delta\in (0,1)$, $p\in[1,\infty]$
and take $\phi$ and $W$ as in Theorem \ref{th_mainconc}. Suppose
that $\bY$ satisfies (GA) and that $t_{\alpha(1-\delta)}$ is a real
number such that $\Proba \paren{ \phi \paren{\moy{\bY}-\mu} >
  t_{\alpha(1-\delta)} }\leq \alpha(1-\delta)$. Then with probability
at least $1-\alpha$, $\phi \paren{\moy{\bY}-\mu}$ is upper bounded by
the minimum between $t_{\alpha(1-\delta)}$ and
\begin{equation} 
\frac{ \E \left[ \phi
      \paren{\moyW{\bY}{W-\moy{W}}} \big| \bY \right] }{B_W} +
  \frac{\norm{\sigma}_{p}}{\sqrt{n}} { \overline{\Phi}}^{-1}
  \paren{\frac{\alpha(1-\delta)}{2}} +
  \frac{\norm{\sigma}_{p}C_W}{nB_W}{\overline{\Phi}}^{-1}
  \paren{\frac{\alpha\delta}{2}}. \label{Equ_concfinalgauss1_comb} 
\end{equation} 
\end{corollary}

\begin{remark}\begin{enumerate}
\item Corollary \ref{cor_mainconc} is a consequence of the proof
  of Theorem \ref{th_mainconc}, rather than of the theorem itself. The point
  here is that $\E \left[ \phi \paren{\moyW{\bY}{W-\moy{W}}}\big|
    \bY\right]$ is almost deterministic, because it concentrates at
  the rate $n^{-1}$ ($=o(n^{-1/2})$). 
\item For instance, if $\phi=\sup(\cdot)$ (resp. $\sup\absj{\cdot}$), Corollary
  \ref{cor_mainconc} may be applied with $t_{\alpha}$ equal to the
  classical Bonferroni threshold for multiple testing (obtained using a simple union bound
  over coordinates)
\begin{equation} \label{eq:seuil_bonf} 
t_{\textrm{Bonf},\alpha} \defeq \frac{1}{\sqrt{n}} \norm{\sigma}_{\infty} \overline{\Phi}^{-1} \paren{ \frac{\alpha}{K}} \bigg( \mbox{resp. } t'_{\textrm{Bonf},\alpha}\defeq \frac{1}{\sqrt{n}} \norm{\sigma}_{\infty} \overline{\Phi}^{-1} \paren{ \frac{\alpha}{2K}}\bigg)
\enspace .\nonumber\end{equation}
 We thus obtain a confidence region almost equal
  to Bonferroni's for small correlations and better than
  Bonferroni's for strong correlations (see simulations in Section
  \ref{sec:simus}). 
\end{enumerate}\end{remark}

The proof of Theorem \ref{th_mainconc} involves results which are of
self interest: the comparison between the expectations of the two
processes  $\E \left[ \phi\paren{\moyW{\bY}{W-\moy{W}}} \big| \bY
\right] $ and $\phi \paren{\moy{\bY} - \mu}$ and the concentration of
these processes around their means. This is examinated in the two
following subsections. The last subsection gives some elements for a
wise choice of resampling weight vectors among several classical
examples. 

\subsection{Comparison in expectation}
In this section, we compare $\E \phi\paren{\moyW{\bY}{W-\moy{W}}}$ and
$\E \phi \paren{\moy{\bY}-\mu}$. We note that these expectations exist
in the Gaussian and the bounded case provided that $\phi$ is
measurable and bounded by a $p$-norm. Otherwise, in particular in
Propositions \ref{prop_espgauss} and \ref{prop_espsym}, we assume that
these expectations exist. In the Gaussian case, these quantities are
equal up to a factor that depends only on the distribution of $W$:
   
\begin{proposition}\label{lemme_espgauss} \label{prop_espgauss} 
Let $\bY$ be a sample satisfying (GA) and $W$ a resampling weight vector.  
Then, for any measurable positive-homogeneous function
$\phi:\R^{K}\rightarrow \R$, we have the following equality  
\begin{equation} 
B_W \E \phi \paren{\moy{\bY}-\mu} = \E \phi\paren{\moyW{\bY}{W-\moy{W}}} \, .
\label{Equ_espgauss} 
\end{equation} 
\end{proposition} 

\begin{remark} 
\begin{enumerate} 
\item In general, we can compute the value of $B_W$ by simulation. For some classical weights, we give bounds or exact expressions in Tab.~\ref{tab:poids_ctes}. 
\item In a non-Gaussian framework, the constant $B_W$ is still relevant, at least asymptotically: in their Theorem 3.6.13, Van der Vaart and Wellner \cite{vdV_Wel:1996} use the limit of $B_W$ when $n$ goes to infinity as a normalizing constant. 
\end{enumerate} 
\end{remark} 
 
When the sample is only symmetric we obtain the following inequalities : 
\begin{proposition}\label{prop_espsym} 
Let $\bY$ be a sample satisfying (SA), $W$ a resampling weight vector and $\phi:\R^{K}\rightarrow \R$ any subadditive, positive-homogeneous function.
\begin{itemize} 
\item[(i)] We have the general following lower bound : 
\begin{equation} 
A_W \E \phi \paren{\moy{\bY}-\mu } \leq \E \phi \paren{\moyW{\bY}{W-\moy{W}}} \, .\label{Equ_majespsym} 
\end{equation} 
\item[(ii)] Moreover, if the weights satisfy the assumption of \eqref{eq:DW}, we have the following upper bound 
\begin{equation} 
D_W \E \phi \paren{\moy{\bY}-\mu} \geq \E \phi \paren{\moyW{\bY}{W-\moy{W}}}.\label{Equ_minespsym} 
\end{equation} 
\end{itemize} 
\end{proposition} 

\begin{remark} 
\begin{enumerate}
\item The bounds \eqref{Equ_majespsym} and \eqref{Equ_minespsym} are
  tight for Rademacher and Random hold-out ($n/2$) weights, but far
  less optimal in some other cases like Leave-one-out (see Section
  \ref{section:weightvectors}).  
\item When $\bY$ is not assumed to be symmetric and $\moy{W}=1$ a.s.,
  Proposition 2 in \cite{Fro:2004} shows that \eqref{Equ_majespsym}
  holds with $\E(W_{1}-\overline{W})_{+}$ instead of $A_W$. Therefore,
  the symmetry of the sample allows us to get a tighter result (for
  instance twice sharper with Efron or Random hold-out ($q$)
  weights). 
\end{enumerate}
\end{remark} 
 
\subsection{Concentration around the expectations} 
In this section we present concentration results for the two processes
$\phi \paren{\moy{\bY}-\mu}$ and $\E \left[
  \phi\paren{\moyW{\bY}{W-\moy{W}}} \big| \bY \right]$ in the Gaussian
framework.  
  
\begin{proposition}\label{prop_concgauss} 
Let $p \in [1,+\infty]$, $\bY$ a sample satisfying (GA) and let $\sigma$ be the vector
$[\var^{1/2}(\bY^{1}_{k})]_k$. Let $\phi: \R^{K}\rightarrow \R$ be any subadditive function, bounded by the $p$-norm.
\begin{itemize} 
\item[(i)] For all $\alpha\in (0,1)$, with probabilty at least
  $1-\alpha$ the following holds: 
\begin{equation} 
\phi \paren{\moy{\bY}-\mu} < \E\phi \paren{\moy{\bY}-\mu} +
\frac{\norm{\sigma}_{p}{\overline{\Phi}}^{-1}(\alpha/2)}{\sqrt{n}} 
\, ,\label{Equ_conc1} 
\end{equation} 
and the same bound holds for the corresponding lower deviations. 
\item[(ii)] Let $W$ be some exchangeable resampling weight vector. 
Then, for all $\alpha\in (0,1)$,  with probabilty at least
  $1-\alpha$ the following holds: 
\begin{equation} 
\E \left[ \phi \paren{ \moyW{\bY}{W-\moy{W}}} \big| \bY \right] <
\E\phi \paren{ \moyW{\bY}{W-\moy{W}}} + \frac{ \norm{\sigma}_{p} C_W
  {\overline{\Phi}}^{-1}(\alpha/2)}{n} 
,\label{Equ_concres1} 
\end{equation} 
and the same bound holds for the corresponding lower deviations. 
\end{itemize} 
\end{proposition} 

The first bound \eqref{Equ_conc1} with a
remainder in $n^{-1/2}$ is classical. The last one
\eqref{Equ_concres1} is much more
interesting since it enlights one of the key properties of the
resampling idea: the ``stabilization''. Indeed, the resampling quantity $\E\left[ \phi
  \paren{\moyW{\bY}{W-\moy{W}} }| \bY \right]$ concentrates around its
expectation at the rate $C_W n^{-1} = o\paren{n^{-1/2}} $ for most of
the weights
(see Section \ref{sec:conc_poids} and Tab.~\ref{tab:poids_ctes} for more details). Thus, compared to the original process, it is almost deterministic and equal to
$B_W\E \phi\paren{\moy{\bY}-\mu}$.

\begin{remark} \label{rem_proconc}
Combining expression \eqref{Equ_espgauss} and Proposition \ref{prop_concgauss} (ii), we derive that for a Gaussian sample $\bY$ and any $p\in[1,\infty]$, the following
  upper bound holds with probability at least $1-\alpha$ :
\begin{equation} 
  \E \norm{\moy{\bY}-\mu}_{p} <  \frac{\E \left[
      \norm{\moyW{\bY}{W-\moy{W}}}_{p} \Big| \bY \right] } {B_W} +
  \frac{\norm{\sigma}_{p} C_W}{n B_W}
  {\overline{\Phi}}^{-1}(\alpha/2) \enspace ,\label{Equ_estrisque}  
\end{equation} 
and a similar lower bound holds.
This gives a control with high probability of the $L^{p}$-risk of the
estimator $\moy{\bY}$ of the mean $\mu\in \R^{K}$ at the rate $C_W B_W^{-1}
n^{-1} $. 
\end{remark} 

\subsection{Resampling weight vectors}  \label{section:weightvectors} \label{sec:conc_poids}  

In this section, we consider the question of choosing some appropriate resampling weight vector $W$ when using Theorem \ref{th_mainconc} or Corollary \ref{cor_mainconc}. We define the following classical resampling weight vectors:
\begin{enumerate} 
\item \textbf{Rademacher}: $W_i$ i.i.d. Rademacher variables, \latin{i.e.} $W_i \in \{-1,1\}$ with equal probabilities.
\item \textbf{Efron}: $W$ has a multinomial distribution with parameters $(n;n^{-1},\ldots, n^{-1})$.
\item \textbf{Random hold-out ($q$)} (R. h.-o.), $q \in \{1, \ldots, n \}$: 
$W_i = \frac{n}{q} \1_{i \in I}$, where $I$ is uniformly distributed on subsets of $\{1,
\ldots, n\}$ of cardinality $q$. These weights may also be called
cross validation weights, or leave-$(n-q)$-out weights. A classical
choice is $q = n/2$ (when $2|n$). When $q=n-1$, these weights are called
\textbf{leave-one-out} weights. 
\end{enumerate} 

\begin{table} \vspace{-5mm}
\begin{center}
\begin{tabular}{|c||c|} 
\hline 
Efron & $2 \paren{1- \frac{1}{n}}^n = A_W \leq B_W \leq \sqrt{\frac{n-1}{n}}$ \quad $C_W = 1$  
\\ Efr., $n \rightarrow + \infty$ & $\frac{2}{e} \leq A_W \leq B_W \leq 1 = C_W$  
\\ \hline  
Rademacher & $1-\frac{1}{\sqrt{n}} \leq A_W \leq B_W \leq \sqrt{ 1-\frac{1}{n}}$ \quad $C_W=1$ \quad $D_W \leq 1+\frac{1}{\sqrt{n}}$ \\  
Rad., $n\rightarrow + \infty$ & $A_W=B_W=C_W=D_W=1$ 
\\ \hline R. h.-o. ($q$) & $A_W = 2 \left(1 - \frac{q}{n}
\right)$ \quad $B_W = \sqrt{ \frac{n}{q} -1 }$ 
\\ 
R. h.-o. ($q$) & $C_W = \sqrt{
  \frac{n}{ n-1}} \sqrt{\frac{n}{q} - 1}$  \quad $D_W = \frac{n}{2q} + \absj{1 - \frac{n}{2q}}$ \\  
R. h.-o. ($n/2$) ($2|n$) & $A_W=B_W=D_W =1$ \quad $C_W=\sqrt{\frac{n}{n-1}}$  \\  
  Leave-one-out & $\frac{2}{n} =A_W \leq B_W =\frac{1}{\sqrt{n-1}}$ \quad$C_W =
\frac{\sqrt{n}}{n-1}$\quad$D_{W}=1$  \\
\hline  
\end{tabular} 
\caption{\label{tab:poids_ctes} Resampling constants for classical
  resampling weight vector.}
\end{center} \vspace{-5mm}
\end{table} 

For these classical weights, exact or approximate values for the quantities $A_W$, $B_W$, $C_W$ and $D_W$ (defined by equations \eqref{eq:AW} to \eqref{eq:DW}) can be easily derived (see Tab. \ref{tab:poids_ctes}).
However, an exact computation of the resampling estimates $\E \left[ \phi \paren{ \moyW{\bY}{W-\moy{W}}} \big| \bY \right]$ using these
weights would be time-consuming when $n$ is large. The more standard way to solve this problem is to compute resampling quantities by Monte-Carlo simulations, \latin{i.e.} picking up a small number of weight vectors (see \cite{Hal:1992}, appendix II for a discussion). But we did not yet investigate the analysis of the corresponding thresholds.

Another way to solve this computation time problem
 is to consider a regular partition $(B_j)_{1 \leq j \leq V}$ 
of $\{1, \ldots, n\}$ (where $V \in \{2, \ldots, n\}$ and $V|n$), and
to define the weights $W_i = \frac{V}{V-1} \1_{i \notin B_J}$ with $J$
uniformly distributed on $\{1, \ldots, V\}$. These weights are called the
\textbf{(regular) $V$-fold cross validation} weights ($V$-f. c.v.),
which are no longer exchangeable but still ``piece-wise
exchangeable''. 
Considering the process
$(\wt{\bY}^{j})_{1 \leq j \leq K}$ where $\wt{\bY}^{j}=\frac{V}{n}\sum_{i\in
B_j}{\bY}^{i}$ is the empirical mean of $\bY$ on block $B_j$, we can
show that Theorem \ref{th_mainconc} can be extended to (regular) $V$-fold
cross validation weights with the following resampling constants:
\[ A_W = \frac{2}{V}\,; \quad 
B_W =\frac{1}{\sqrt{V-1}}\,; \quad 
C_W = \sqrt{n}(V-1)^{-1}\,; \quad D_W = 1 \enspace . \]
When $V$ does not divide $n$ and the blocks are no longer regular, Theorem \ref{th_mainconc}
can also be generalized, but the constants have more complex
expressions.

Note that in the Gaussian 
framework of \eqref{Equ_estrisque}, 
$V$-fold cross-validation weights approximate the estimation risk
$\E \norm{\moy{\bY}-\mu}_p$ by $\frac{\sqrt{V-1}}{V^2} \sum_{j=1}^V
\norm{ \wt{\bY}^{(-j)} - \wt{\bY}^{j} }_p $, where $\wt{\bY}^{(-j)}$
is the mean of the $(\wt{\bY}^{\ell})_{\ell\neq j}$\,; which bears
a strong analogy with the usual cross-validation philosophy.
Actually, the ``classical'' leave-one-out estimator
$\frac{1}{n}\sum_{i=1}^n \norm{ \wt{\bY}^{(-i)} - {\bY}^{i} }_p $ 
approximates a different quantity, the prediction risk $\E\norm{\moy{\bY}-\bY^{n+1}}_p$
for a new independent vector $\bY^{n+1}$\,. However, under (GA)
the two types of risk are proportional, $\sqrt{n+1} \E
\norm{\moy{\bY}-\mu}_p=\E\norm{\moy{\bY}-\bY^{n+1}}_p$\,;
taking into account this scaling we conclude that
our estimator (with $V=n$) coincides with the classical leave-one-out 
(up to the factor $\sqrt{1-1/n^2}\sim 1$).
To guide our choice for a specific resampling scheme, 
the first comparison point is that $t_{\alpha,W}(\bY)$ should be
an accurate upper bound of the ideal threshold. 
Under the Gaussian
assumption, in view of \eqref{Equ_concfinalgauss1}, $C_W B_W^{-1}$
appears as a relevant accuracy index for $t_{\alpha,W}$. However,
a second comparison point is the price of
an exact computation of $t_{\alpha,W}$ in practice. Since one
must consider each possible weight vector to compute exactly the
threshold, we use the cardinality of the support of $\mathcal{L}(W)$
as a complexity index.  

As shown in Tab.~\ref{tab:poids_ctes_2}, there is an
\emph{accuracy-complexity trade-off} for choosing the weights. 
Since for all exchangeable weights $C_W B_W^{-1} \geq \sqrt{n/(n-1)}$, R. h.-o.($n/2$) and leave-one-out weights are optimal for accuracy (Rademacher and Efron being "almost optimal"). 
On the other hand, $V$-fold c.-v. is less accurate, losing a factor  $\sqrt{(n-1)/(V-1)}$.
On the computational viewpoint, the leave-one-out is the only reasonable
exchangeable procedure (at least when $n$ and $K$ are large), and 
$V$-f. c.v. looks 
even  more attractive. Considering that
$t_{\alpha,W}$ involves the sum of terms of order $C_W B_W^{-1} n^{-1}$ and
$n^{-1/2}$, the best choice of $V$ should be rather small for most
applications. We do not give here any universal optimal $V$ since it
does not exist, but we suggest to use Tab.~\ref{tab:poids_ctes_2} to
choose it.

\begin{table} \vspace{-5mm}
\begin{center} 
\begin{tabular}{|c||l|c|} 
\hline Resampling & $\:\:\:\:\:\:C_W B_W^{-1}$ (accuracy) & $\card \paren{\supp \mathcal{L}(W)}$ (complexity) \\ \hline \hline  
Efron & $\leq \frac{1}{2} \paren{1-\frac{1}{n}}^{-n} \xrightarrow[n \rightarrow \infty]{} \frac{e}{2} $ & $n^n$ 
\\ \hline 
 Rademacher & $\leq \paren{1 - n^{-1/2}}^{-1} \xrightarrow[n \rightarrow \infty]{} 1$ & $2^n$  
\\ \hline R. h.-o. ($n/2$) & $= \sqrt{ \frac{n}{ n-1}} \xrightarrow[n \rightarrow \infty]{} 1$ & $\binom{n}{n/2} \propto n^{-1/2} 2^n$ \\ 
Leave-one-out & $= \sqrt{ \frac{n}{ n-1}} \xrightarrow[n \rightarrow \infty]{} 1$ & $n$  
\\ \hline \hline  
regular $V$-fold c.-v. & $= \sqrt{\frac{n}{V-1}} $  & $V$ 
\\ \hline  
\end{tabular} 
\caption{\label{tab:poids_ctes_2} Choice of the resampling weight vectors : accuracy-complexity tradeoff.} 
\end{center} \vspace{-7mm}
\end{table} 

\section{Confidence region using resampled quantiles}

In the previous section we have shown how to derive non-asymptotic
confidence regions for the mean of a Gaussian (resp. bounded) vector with unknown
correlation structure; for this we used a concentration property
of the quantities $\phi(\ol{\bY}-\mu)$ and $\E \left[
  \phi(\ol{\bY}_{[W-\ol{W}]}) | \bY \right]$ around their mean. The
Gaussian  (resp. McDiarmid's)
concentration property allowed us to bound deviations from
this mean by the deviations of a suitably scaled normal
(resp. subgaussian) variable.
Through this approach, the level of the
confidence region is rigorously controlled for any fixed sample size.

However, the obtained confidence regions are somewhat unsatisfying
because they appear to be too conservative in practice. The principal
reason for this is that $\phi(\ol{\bY}-\mu)$ is of
course not 
a Gaussian variable (even when $\bY$ is)
. Therefore, in spite
of the power of the Gaussian concentration property, using Gaussian
tails as a bound for the deviations of the above non-Gaussian variable must necessarily
result in losing some slack.

On the other hand, in most applications of resampling procedures, it is common to estimate
the quantiles of a variable like $\phi\paren{\moy{\bY}-\mu}$ by the
quantiles of the corresponding resampled distribution
$\mathcal{L} \paren{ \phi\paren{\moyW{\bY}{W - \moy{W}}} \big| \bY }$,
and to use these quantiles to construct 
a confidence region. Again, while many asymptotic results are
available to justify this method (for instance \cite{vdV_Wel:1996}), our goal here is to
derive a non-asymptotic region based on a similar approach
for which the confidence level
is proved to hold for any fixed sample size.

For this we apply a principle that is close in spirit to
exact tests, \latin{i.e.} by taking advantage of an invariance property
(here symmetry around the mean) of the initial distribution and
using a resampling scheme that respects this invariance.
For this reason the scope of the current section is
far less general: instead of
covering generic resampling weights, we only consider 
the particular Rademacher resampling scheme. Let us define for a function $\phi$ the resampled empirical quantile:
\[
q_{\alpha}(\phi,\bY) = 
\inf\set{x\in\mbr \telque \probb{W}{\phi(\ol{\bY}_{[W]})>x}\leq
  \alpha}\,, 
\]
wherein  $W$  is an i.i.d Rademacher weight vector. 
We now state the main technical result of this section:

\begin{proposition}
\label{le:quantiles}\label{prop:quantiles}
Fix $\delta,\al \in
(0,1)$. Let $\bY$ be a data sample satisfying assumption (SA).
Let $f: \paren{\mbr^K}^n  \rightarrow [0,\infty )$ be a nonnegative
(measurable) function on the set of data samples. 
Let $\phi$ be a {\em nonnegative}, subadditive, positive-homogeneous function.
Denote $\wt{\phi}(x) = \max\paren{\phi(x),\phi(-x)}$\,. Finally, for $\eta \in (0,1)$\,, denote
\[
\ol{\cB}(n,\eta) = \min\set{k \in \{0,\dots,n\} \telque 2^{-n} \sum_{i=k+1}^n \binom{n}{i} < \eta}\,,
\]
the upper quantile function of a binomial $(n,\frac{1}{2})$ variable.
Then we have:
\begin{multline}
\prob{\phi(\overline{\bY}-\mu)>q_{\alpha(1-\delta)} \paren{\phi,\bY - \moy{\bY} }+f(\bY)}\\
\leq \alpha + \prob{\wt{\phi}(\overline{\bY}-\mu)> \frac{n}{2 \ol{\cB}\paren{n,\frac{\al\delta}{2}}-n}f(\bY)} \nonumber
\end{multline}
\end{proposition}

\begin{remark} By Hoeffding's inequality, $\frac{n}{2 \ol{\cB}\paren{n,\frac{\al\delta}{2}}-n} \geq \paren{\frac{n}{2\ln\paren{\frac{2}{\al\delta}}}}^{1/2}.$ 
\end{remark}
By iteration of this proposition we obtain the following corollary:

\begin{corollary}
\label{cor:quantregiongauss}
Fix $J$ a positive integer, $(\alpha_i)_{i=0,\ldots,J-1}$ a finite sequence in $(0,1)$ and $\beta,\delta \in (0,1)$\,. 
Let $\bY$ be a data sample satisfying assumption (SA). Let
$\phi:\R^{K}\rightarrow \R$ be a nonnegative, subadditive, positive-homogeneous function
and $f: \paren{\mbr^K}^n  \rightarrow [0,\infty )$ be a nonnegative
function on the set of data samples. 
Then the following holds:
\begin{multline}
\label{controle_quant}
\prob{\phi(\ol{\bY}-\mu) >  q_{(1-\delta)\alpha_0}(\phi,\bY-\ol{\bY})
  + \sum_{i=1}^{J-1} \gamma_i
  {q}_{(1-\delta)\alpha_i}(\wt{\phi},\bY-\ol{\bY}) 
  +\gamma_{J} f(\bY) }
\\ \leq \sum_{i=0}^{J-1} \alpha_i +
\prob{\wt{\phi}(\ol{\bY}-\mu) > f(\bY)} \,,
\end{multline}
where, for $k\geq 1$, $\displaystyle\gamma_k = 
n^{-k}\prod_{i=0}^{k-1} \paren{2\ol{\cB}\paren{n,\frac{\al_i\delta}{2}}-n}\,.
$
\end{corollary}
The rationale behind this result is that the sum appearing 
inside the probability should be interpreted as a series of
corrective terms of decreasing order
of magnitude, since we expect the sequence $\gamma_k$
to be sharply decreasing. Looking at Hoeffding's bound, this will be the case
if the levels are such that $\alpha_i \gg \exp (-n)$\,. 

Then comes the remaining issue of the trailing term on the
right-hand-side. While it is tempting to think that it would be
possible to obtain a self-contained result based on the symmetry assumption
(SA) alone, we did not succeed in this direction. To upper-bound the trailing term,
we can assume some additional regularity assumption on the distribution of the data.
For example, if the data are Gaussian or bounded, we can apply the results in the previous section 
(or apply some other device like Bonferroni's bound \eqref{eq:seuil_bonf}).
The point is that this bound does not have to be particularly sharp, since
we expect (in favorable cases) the trailing probability term on the right-hand side
as well as the contribution of $\gamma_J f(\bY)$ to the left-hand side
to be almost negligible.

It seems
plausible that at least a minor regularity assumption 
(supposedly significantly weaker than assuming a Gaussian distribution or 
bounded data) is actually a
necessary condition in addition to (SA) to obtain a self-contained bound and
ensure that nothing pathological happens with the
extreme quantiles, but this remains as an interesting open issue.

As before, for computational reasons, it might be relevant to consider a block-wise Rademacher resampling
scheme.
\section{Simulations} \label{sec:simus}
For simulations we consider data of the form $Y_t = \mu_t + G_t$\,,
where $t$ belongs to an $m\times m$ discretized 2D torus of $K=m^2$ ``pixels'',
identified with $\mbt^2_m=(\mbz/m\mbz)^2$\,, 
and $G$ is a centered Gaussian vector obtained by 2D discrete
convolution of an i.i.d. standard Gaussian field 
(``white noise'') on $\mbt^2_m$ with a function $F:\mbt^2_m \rightarrow \mbr $ 
such that $\sum_{t\in\mbt^2_m}F^2(t)=1$
\,. This ensures that $G$
is a stationary Gaussian process on the discrete torus, it is
in particular isotropic with $\e{G_t^2}=1$ for all $t\in\mbt^2_m$\,.

In the simulations
below we consider for the function $F$ a ``Gaussian''
convolution filter of bandwith $b$ on the torus:
\[
F_b(t)  = C_b \exp\paren{-d(0,t)^2/b^2}\,,
\]
where $d(t,t')$ is the standard distance on the torus and $C_b$ is a normalizing
constant. Note that for actual
simulations it is more convenient to work in the Fourier domain and 
to apply the inverse DFT which can be computed efficiently.
We then compare
the different thresholds obtained by the methods proposed in this work
for varying values of $b$\,. Remember that the only information available to
the algorithm is the bound on the marginal variance; the form of the function $F_b$
itself is of course unknown. 

\begin{figure}[t]
\label{restexture}
\raisebox{4.5mm}{\includegraphics[height=40mm]{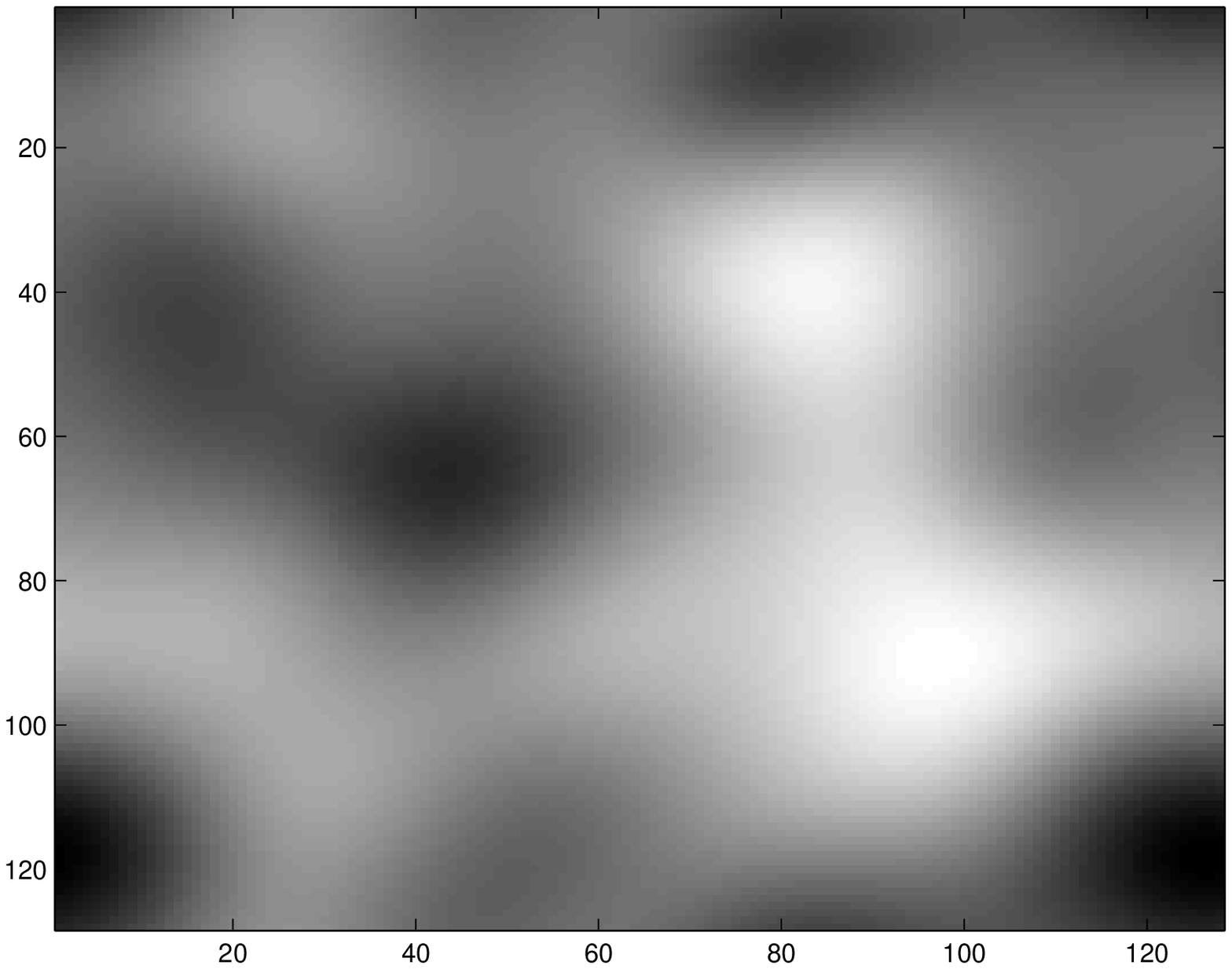}}
\hfill \includegraphics[height=50mm]{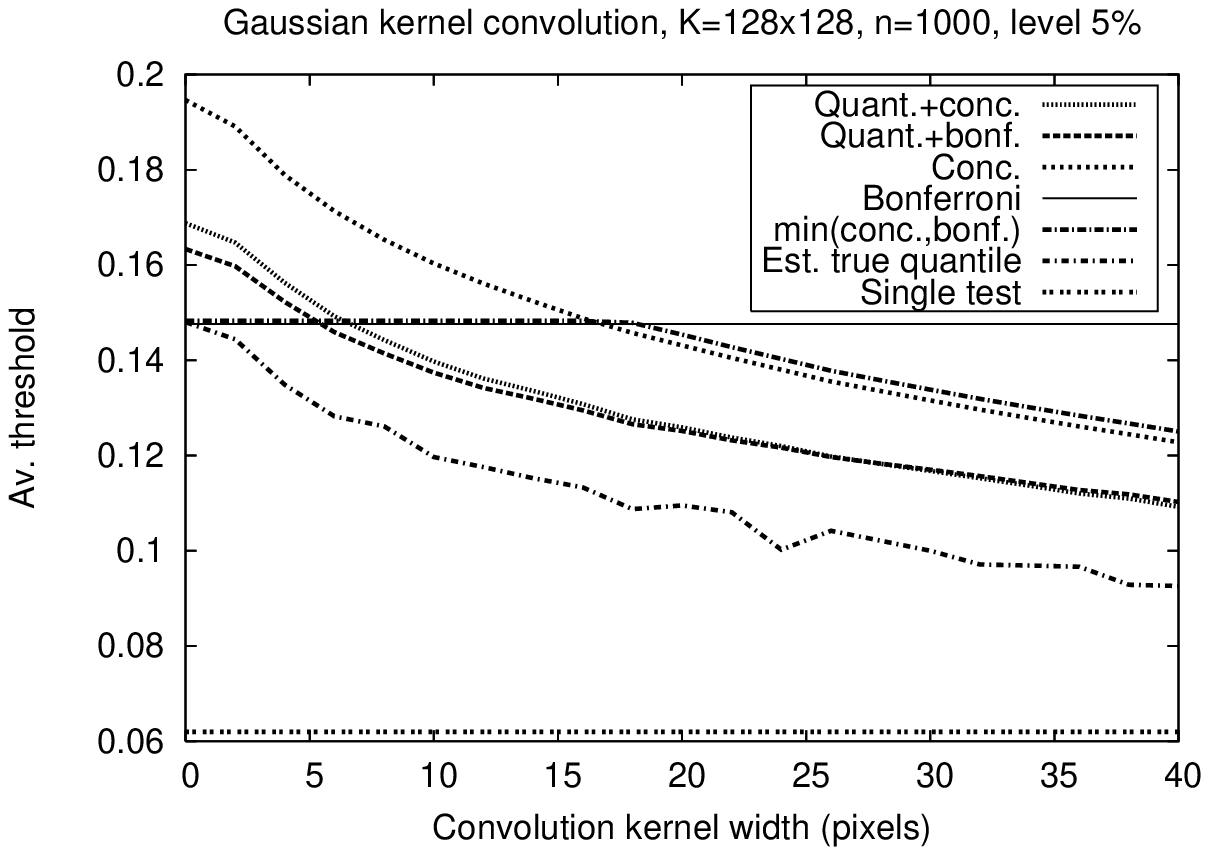}
\caption{Left: example of a 128x128 pixel image obtained by convolution of
Gaussian white noise with a (toroidal) Gaussian filter with width $b=18$ pixels.
Right: average thresholds obtained for the different approaches, see text.}
\vspace{-4mm}
\end{figure}

On Fig.~\ref{restexture} we compare the thresholds obtained when
$\phi = \sup \absj{\cdot}$\,, which corresponds to the two-sided multiple testing
situation.
We use the different approaches proposed in this work, with the following
parameters: the dimension is $K=128^2=16384$\,, the number of data points per sample
is $n=1000$ (much smaller than $K$, so that we really are in a non-asymptotic framework), the width $b$ takes even values in the range $[0,40]$\,,
the overall level is $\alpha = 0.05$\,. For the concentration threshold \eqref{Equ_concfinalgauss1} ('conc.'), 
we used Rademacher weights. For the ``compound'' threshold of Corollary \ref{cor_mainconc} ('min(conc.,bonf.)'), we used $\delta = 0.1$ and the Bonferroni threshold
$t'_{\textrm{Bonf},0.9\alpha}$ as the deterministic reference threshold. For the quantile approach 
\eqref{controle_quant}, 
we used $J=1$\,, $\alpha_0 = 0.9\alpha$\,, $\delta = 0.1$\,, and the function
$f$ is given either by the Bonferroni threshold ('quant.+bonf.') or the concentration
threshold ('quant.+conc.'), both at level $ 0.1\alpha$\,.  Each point 
represents an average over 50 experiments. Finally, we included in
the figure the Bonferroni threshold $t'_{\textrm{Bonf},\alpha}$
, the threshold for a single test for
comparison, and an estimation of the true quantile (actually, an
empirical quantile over 1000 samples).

The quantiles or expectation with Rademacher weights were estimated by 
Monte-Carlo with 1000 draws. On the figure we did not include standard
deviations: they are quite low, of the order of $10^{-3}$\,, although
it is worth noting that the quantile threshold has a standard deviation roughly
twice as large as the concentration threshold (we did not investigate
at this point what part
of this variation is due to the MC approximation).

The overall conclusion of this preliminary
experiment is that the different
thresholds proposed in this work are relevant in the sense that they
are smaller than the Bonferroni threshold provided the vector has
strong enough correlations. As expected,
the quantile approach appears to lead to tighter thresholds. (However, this might
not be always the case for smaller sample sizes.) One advantage of the concentration
approach is that the 'compound' threshold \eqref{Equ_concfinalgauss1_comb}
can ``fall back'' on the Bonferroni threshold when needed, at the price
of a minimal threshold increase. 
\section{Proofs}
 
\begin{proof}[Proof of Proposition \ref{prop_espgauss}] 
Denoting by $\mathbf{\Sigma}$ the common covariance matrix of the $\bY^i$, we have $\mathcal{L}(\moyW{\bY}{W-\moy{W}}|W) = (n^{-1} \sum_{i=1}^{n}(W_{i}-\moy{W})^{2})^{1/2} \mathcal{N}(0,n^{-1}\mathbf{\Sigma})$, and the result follows because $\mathcal{L}(\moy{\bY}-\mu)=\mathcal{N}(0,n^{-1}\mathbf{\Sigma})$ and $\phi$ is positive-homogeneous. \qed
\end{proof}

\begin{proof}[Proof of Proposition \ref{prop_espsym}] (i). 
By independence between $W$ and $\bY$, using the positive homogeneity, then convexity of $\phi$, 
for every realization of $\bY$ we have: 
\begin{align*} 
A_W \phi\paren{\moy{\bY}-\mu} &= \phi \paren{ \E \left[ \frac{1}{n} \sum_{i=1}^{n} \absj{W_{i}-\moy{W}}  \paren{\bY^{i}-\mu} \bigg| \bY \right] }\\ 
&\leq  \E \left[ \phi \paren{\frac{1}{n} \sum_{i=1}^{n} \absj{W_{i}-\moy{W}} \paren{\bY^{i}-\mu} } \bigg| \bY \right] \, . 
\end{align*} 

We integrate with respect to $\bY$, and use the symmetry of the $\bY^{i}$ with respect to $\mu$ and again the independence between $W$ and $\bY$ to show finally that 
\begin{align*} 
A_W \E \left[ \phi\paren{\moy{\bY}-\mu} \right] &\leq 
\E \left[ \phi \paren{ \frac{1}{n} \sum_{i=1}^{n} \absj{W_{i}-\moy{W}} \paren{ \bY^{i}-\mu } } \right] \\ 
&=\E \left[ \phi \paren{ \frac{1}{n} \sum_{i=1}^{n} \paren{ W_{i}-\moy{W} } \paren{ \bY^{i}-\mu } } \right] = \E \left[ \phi \paren{ \moyW{\bY}{W-\moy{W}} } \right] \, . 
\end{align*} 
We obtain (ii) via the triangle inequality and the same symmetrization trick. \qed
\end{proof} 
 
\begin{proof}[Proof of Proposition \ref{prop_concgauss}] 
We denote by $\bA$ a square root of the common covariance matrix of
the $\bY^{i}$ and by $(a_{k})_{1 \leq k\leq K}$ the rows of $\bA$. 
If $\bG$ is a $K\times m$ matrix with standard centered i.i.d. Gaussian entries,
then $\bA \bG$ has the same distribution as $\bY-\mu$\,. 
We
let for all $\zeta \in \paren{{\R}^{K}}^{n}$, $T_{1}(\zeta) \defeq
\phi \paren{\frac{1}{n}\sum_{i=1}^{n} \bA\zeta_{i} }$ and
$T_{2}(\zeta)\defeq \E \phi\paren{\frac{1}{n}\sum_{i=1}^{n}
  (W_{i}-\moy{W}) \bA\zeta_{i} }$.  
From the Gaussian concentration theorem of Cirel'son, Ibragimov and Sudakov
(see for example \cite{Mas:2003:St-Flour}, Theorem 3.8), we just need
to prove that $T_{1}$ (resp. $T_{2}$) is a Lipschitz function with
constant $\norm{\sigma}_{p}/\sqrt{n}$
(resp. $\norm{\sigma}_{p}C_W/{n}$), for the Euclidean norm $\norm{\cdot}_{2,Kn}$ on
$\paren{{\R}^{K}}^{n}$.
Let $\zeta,\zeta' \in \paren{{\R}^{K}}^{n}$.
Using Cauchy-Schwartz's inequality coordinate-wise and
$\norm{a_{k}}_{2}=\sigma_{k}$, we deduce  
\begin{equation*} 
\absj{T_{1}(\zeta)-T_{1}(\zeta')}\leq \norm{ \frac{1}{n}\sum_{i=1}^{n} \bA \paren{\zeta_{i}-\zeta'_{i}} }_{p}\leq \norm{\sigma}_{p} \norm{\frac{1}{n}\sum_{i=1}^{n} \paren{\zeta_{i}-\zeta'_{i}}}_{2}.
\end{equation*} 
Therefore, we get
$\absj{T_{1}(\zeta)-T_{1}(\zeta^{\prime})}\leq \frac{\norm{\sigma}_{p}}{\sqrt{n}} \norm{\zeta-\zeta^{\prime}}_{2,Kn}$  
 by convexity of $x\in{\R}^{K}\rightarrow \norm{x}_{2}^{2}$, and we obtain (i). For $T_{2}$, we use the same method as for $T_{1}$ : 
\begin{align} 
\absj{T_{2}(\zeta)-T_{2}(\zeta')}&\leq \norm{\sigma}_{p}  \E \norm{ \frac{1}{n}\sum_{i=1}^{n} (W_{i}-\moy{W})(\zeta_{i}-\zeta'_{i})}_{2} \notag \\ 
&\leq \frac{\norm{\sigma}_{p}}{n}  \sqrt{\E \norm{\sum_{i=1}^{n} (W_{i}-\moy{W})(\zeta_{i}-\zeta'_{i})}_{2}^{2}}\label{equ_lipT2_1} \enspace . 
\end{align} 
We now develop $\norm{\sum_{i=1}^{n}
  (W_{i}-\moy{W})(\zeta_{i}-\zeta'_{i})}_{2}^{2}$ in the Euclidean
space $\R^K$ (note that from $\paren{\sum_{i=1}^{n}
  (W_{i}-\moy{W})}^2=0$, we have
$\E(W_{1}-\moy{W})(W_{2}-\moy{W})=-C_W^{2}/n$) : 
\begin{align*} 
\E \norm{\sum_{i=1}^{n} (W_{i}-\moy{W})(\zeta_{i}-\zeta'_{i})}_{2}^{2} &= C_W^{2}(1-1/n)\sum_{i=1}^n  
\norm{\zeta_i - \zeta^{\prime}_i}_{2}^2  -\frac{C_W^{2}}{n} \sum_{i \neq j} <{\zeta_{i}-\zeta'_{i}},{\zeta_{j}-\zeta'_{j}}> \\ 
&= C_W^{2}\sum_{i=1}^n \norm{\zeta_i - \zeta^{\prime}_i}_{2}^2 -\frac{C_W^{2}}{n} \norm{\sum_{i=1}^{n} (\zeta_{i}-\zeta'_{i})}_{2}^{2}\enspace . 
\end{align*} 
Consequently, 
\begin{align} 
\E \norm{ \sum_{i=1}^{n} \paren{ W_{i}-\moy{W}} \paren{ \zeta_{i}-\zeta'_{i} } }_{2}^{2} \leq C_W^{ 2} \sum_{i=1}^n \norm{\zeta_i - \zeta^{\prime}_i}_{2}^2 
\leq C_W^{2} \norm{\zeta-\zeta'}_{2,Kn}^{2}.\label{equ_lipT2_2} 
\end{align} 
Combining expression (\ref{equ_lipT2_1}) and (\ref{equ_lipT2_2}), we find that $T_{2}$ is $\norm{\sigma}_{p}C_W/{n}$-Lipschitz. \qed
\end{proof} 
 
\begin{proof}[Proof of Theorem \ref{th_mainconc}] 
The case (BA)($p,M$) and (SA) is obtained by combining Proposition \ref{prop_espsym} and McDiarmid's inequality (see for instance \cite{Fro:2004}). 
 The (GA) case is a straightforward consequence of Proposition \ref{prop_espgauss} and the proof of Proposition \ref{prop_concgauss}. \qed
\end{proof} 
 
\begin{proof}[Proof of Corollary \ref{cor_mainconc}] 
  From Proposition \ref{prop_concgauss} (i), with probability at least
  $1-\alpha(1-\delta)$, $\phi \paren{\moy{\bY}-\mu}$ is upper bounded by
  the minimum between $t_{\alpha(1-\delta)}$ and $\E\phi
  \paren{\moy{\bY}-\mu} +
  \frac{\norm{\sigma}_{p}{\overline{\Phi}}^{-1}(\alpha(1-\delta)/2)}{\sqrt{n}}$
  (because these thresholds are deterministic).
  In addition, Proposition \ref{prop_espgauss} and Proposition
  \ref{prop_concgauss} (ii) give that with probability at least
  $1-\alpha\delta$, $\E \phi \paren{\moy{\bY}-\mu} \leq
  \frac{\E\paren{\phi \paren{\moy{\bY}-\mu}|\bY}}{B_W} +
  \frac{\norm{\sigma}_{p}C_{W}}{B_{W}n}{\overline{\Phi}}^{-1}(\alpha\delta/2)$.
The result follows by combining the two last expressions. \qed 
\end{proof}

\begin{proof}[Proof of Proposition \ref{le:quantiles}] Remember the following inequality coming from the definition of the
quantile $q_\alpha$\,: for any fixed $\bY$
\begin{equation}
\label{eq:inegquant}
\probb{\gru}{\phi\paren{\overline{\bY}_{[\gru]}}>q_{\alpha}(\phi,\bY)}
\leq \alpha \leq \probb{\gru}{\phi\paren{\overline{\bY}_{[\gru]}}\geq 
q_{\alpha}(\phi,\bY)}\,,
\end{equation}
which will be useful in this proof. We have 
\begin{align}
\probb{\bY}{\phi(\overline{\bY}-\mu)>q_{\alpha}(\phi, \bY-\mu)}
&=
\ee{\gru}{\probb{\bY}{\phi\big(\overline{(\bY-\mu)}_{[\gru]}\big)>q_{\alpha}(\phi,(\bY-\mu)_{[\gru]})}}  
\nonumber \\
&=\ee{\bY}{\probb{\gru}{\phi\paren{\overline{(\bY-\mu)}_{[\gru]}}>q_{\alpha}(\phi,\bY-\mu)}}
\nonumber \\
&\leq \alpha\,.
\label{true_recent}
\end{align}
The first equality is due to the fact that the distribution of
$\bY$ satisfies assumption (SA), hence the distribution of $(\bY-\mu)$ invariant by reweighting
by (arbitrary) signs $\gru \in \set{-1,1}^n$\,. In the second
equality we used Fubini's theorem and the fact that for any arbitrary
signs $\gru$ as above
$q_{\alpha}(\phi,(\bY-\mu)_{[\gru]}) =  q_{\alpha}(\phi,\bY-\mu)$\,; finally the last
inequality comes from \eqref{eq:inegquant}. Let us define the event
\[
\Om = \set{\bY \telque q_\al(\phi,\bY-\mu) \leq q_{\al(1-\delta)}(\phi,\bY-\ol{\bY}) + f(\bY) }\,;
\]
then we have using \eqref{true_recent} :
\begin{eqnarray}
\prob{\phi(\ol{\bY}-\mu) > q_{\al(1-\delta)}(\phi,\bY-\ol{\bY}) + f(\bY)}&\leq& \prob{ \phi(\ol{\bY}-\mu) > q_\al(\phi,\bY-\mu)} + \prob{\bY \in \Om^c}\nonumber\\
&\leq& \al + \prob{\bY \in \Om^c} \enspace . \label{to_om}
\end{eqnarray}

We now concentrate on the event $\Omega^c$\,. Using the subadditivity
of $\phi$, and the fact that $\ol{(\bY-\mu)}_{[\gru]} = \ol{(\bY-\ol{\bY})}_{[\gru]}
+ \ol{\gru}(\ol{\bY}-\mu)$\,, we have for any fixed $\bY\in\Om^c$:
\begin{align}
\al & \leq \probb{\gru}{\phi(\ol{(\bY-\mu)}_{[\gru]}) \geq q_{\al}(\phi,\bY-\mu)} \nonumber \\
& \leq \probb{\gru}{\phi(\ol{(\bY-\mu)}_{[\gru]}) > q_{\al(1-\delta)}(\phi,\bY-\ol{\bY}) +
  f(\bY)} \nonumber \\
&\leq \probb{\gru}{\phi(\ol{(\bY-\ol{\bY})}_{[\gru]}) >
  q_{\al(1-\delta)}(\phi,\bY-\ol{\bY})} +
\probb{\gru}{\phi(\ol{\gru}(\ol{\bY}-\mu)) > f(\bY)} \nonumber\\
&\leq \al(1-\delta) + \probb{\gru}{\phi(\ol{\gru}(\ol{\bY}-\mu)) >
  f(\bY)}\,.\nonumber
\end{align}
For the first and last inequalities we have used \eqref{eq:inegquant},
and for the second inequality the definition of $\Omega^c$.
From this we deduce that
\[
\Om^c \subset \set{\bY \telque \probb{\gru}{\phi(\ol{\gru}(\ol{\bY}-\mu)) >  f(\bY)} \geq \al\delta}\,.
\]
Now using the homogeneity of $\phi$, and the fact that both
$\phi$ and $f$ are nonnegative:
\begin{align*}
\probb{\gru}{\phi(\ol{\gru}(\ol{\bY}-\mu)) > f(\bY)} 
& = \probb{\gru}{\abs{\ol{\gru}} > \frac{f(\bY)}{{\phi}(\mathrm{sign}(\ol{\gru})(\ol{\bY}-\mu))}}\\
& \leq  \probb{\gru}{\abs{\ol{\gru}} > \frac{f(\bY)}{\wt{\phi}(\ol{\bY}-\mu)}}\\
& = 2 \prob{ \frac{1}{n}(2B_{n,\frac{1}{2}} - n) > \frac{f(\bY)}{\wt{\phi}(\ol{\bY}-\mu)} \bigg| \bY }\,,
\end{align*}
where $B_{n,\frac{1}{2}}$ denotes a binomial $(n,\frac{1}{2})$ variable (independent of $\bY$).
From the two last displays we
conclude
\[
\Om^c \subset \set{\bY \telque \wt{\phi}(\ol{\bY}-\mu) >
\frac{n}{2\ol{\cB}\paren{n,\frac{\al\delta}{2}}-n}
f(\bY)}\,,
\]
which, put back in \eqref{to_om}, leads to the desired conclusion. \qed
\end{proof}

\section*{Acknowledgements}
We want to thank Pascal Massart for his particulary relevant suggestions.
\bibliographystyle{alpha}
\bibliography{ABR07}

\begin{thebibliography}{VdVW96}

\bibitem[Efr79]{Efr:1979}
B.~Efron.
\newblock Bootstrap methods: another look at the jackknife.
\newblock {\em Ann. Statist.}, 7(1):1--26, 1979.

\bibitem[Fro04]{Fro:2004}
Magalie Fromont.
\newblock Model selection by bootstrap penalization for classification.
\newblock In {\em Learning theory}, volume 3120 of {\em Lecture Notes in
  Comput. Sci.}, pages 285--299. Springer, Berlin, 2004.

\bibitem[Hal92]{Hal:1992}
Peter Hall.
\newblock {\em The bootstrap and {E}dgeworth expansion}.
\newblock Springer Series in Statistics. Springer-Verlag, New York, 1992.

\bibitem[Mas05]{Mas:2003:St-Flour}
Pascal Massart.
\newblock Concentration inequalities and model selection (lecture notes of the
  {S}t-{F}lour probability summer school 2003).
\newblock Available online at {\tt{http://www.math.u-psud.fr/\~{
  }massart/stf2003\_massart.pdf}}, 2005.

\bibitem[VdVW96]{vdV_Wel:1996}
Aad~W. Van~der Vaart and Jon~A. Wellner.
\newblock {\em Weak convergence and empirical processes}.
\newblock Springer Series in Statistics. Springer-Verlag, New York, 1996.
\newblock With applications to statistics.

\end{thebibliography}
\end{document}